\documentclass[12pt]{amsart}
\usepackage{amssymb, amsmath, amscd, dcpic, pictexwd}
\usepackage{hyperref}
\usepackage[all]{xy}

\newtheorem{thm}{Theorem}[section]

\newtheorem{cor}[thm]{Corollary}
\newtheorem{lem}[thm]{Lemma}
\newtheorem{prop}[thm]{Proposition}
\newtheorem{defn}[thm]{Definition}
\newtheorem{rem}[thm]{Remark}
\newtheorem{exm}[thm]{Example}
\numberwithin{equation}{section}

\newcommand{\cA}{\mathcal A}

\newcommand{\cD}{\mathcal D}

\newcommand{\cS}{\mathcal S}
\newcommand{\cO}{\mathcal{O}}
\newcommand{\cL}{\mathcal{L}}

\newcommand{\cE}{\mathcal E}

\newcommand{\cM}{\mathcal M}
\newcommand{\cH}{\mathcal H}

\newcommand{\cY}{\mathcal{Y}}

\newcommand{\fg}{\mathfrak{g}}
\newcommand{\fp}{\mathfrak{p}}
\newcommand{\ft}{\mathfrak{t}}
\newcommand{\fh}{\mathfrak{h}}
\newcommand{\fb}{\mathfrak{b}}

\newcommand{\fk}{\mathfrak{k}}

\newcommand{\HH}{\mathbb{H}}

\newcommand{\CC}{\mathbb{C}}
\newcommand{\PP}{\mathbb{P}}
\newcommand{\ZZ}{\mathbb Z}
\newcommand{\into}{\hookrightarrow}
\newcommand{\onto}{\twoheadrightarrow}

\newcommand{\ra}{\rightarrow}

\newcommand{\bra}{{\langle}}
\newcommand{\ket}{{\rangle}}

\newcommand{\C}{\CC}

\newcommand{\RR}{{\mathbb R}}

\newcommand{\Aut}{{\mbox{Aut~}}}
\newcommand{\End}{{\mbox{End~}}}
\newcommand{\rank}{{\mbox{rank}}}

\newcommand{\bs}{\bigskip}
\newcommand{\cut}{\vskip-.3in}

\def\comment#1{}

\author{Bong H. Lian, Ruifang Song, and Shing-Tung Yau}

\title{Period Integrals and Tautological Systems}
\begin{document}
\maketitle
\begin{abstract}
We study period integrals of CY hypersurfaces in a partial flag variety. We construct a regular holonomic system of differential equations which govern the period integrals. By means of representation theory, a set of generators of the system can be described explicitly.
The results are also generalized to CY complete intersections. The construction of these new systems of differential equations have lead us to the notion of a tautological system.

\end{abstract}


\tableofcontents 
\baselineskip=16pt plus 1pt minus 1pt
\parskip=\baselineskip

\section{Introduction}

Let $X$ be a nonsingular $d$-dimensional Fano variety, i.e. $K_X^{-1}$ is ample. Assume that a general section $f_0\in V^*=H^0(X,K_X^{-1})$ defines a nonsingular CY variety $Y_{f_0}=\{f_0=0\}$. The local Torelli theorem implies that the line $H^{d-1,0}(Y_f)\subset H^{d-1}(Y_{f_0})$ determines the isomorphism class of $Y_f$, for $f$ close to a fixed $f_0$. Period integrals provide a way to parameterize the lines $H^{d-1,0}(Y_f)$, as $f$ vary. By the Kodaira-Nakono-Akizuki vanishing theorem and Serre duality, the Poincar\'e residue sequence collapses to give an isomorphism 
$$Res:H^0(X,\Omega^d(Y_f))\ra H^{d-1,0}(Y_f).$$ 
Let $\gamma_i$ be a basis of the free part of $H_{d-1}(Y_{f_0},\ZZ)$, and $\tau(\gamma_i)\in H_d(X,\ZZ)$ be a small tube over $\gamma_i$. We can choose a local family of meromorphic $d$-forms $\Omega_f$ with a single pole along $Y_f$, so that
$$
\int_{\gamma_i} Res~\Omega_f=\int_{\tau(\gamma_i)}\Omega_f.
$$
These period integrals determine the line $H^{d-1,0}(Y_f)$, and they fit together to form a sheaf of functions, which we call the period sheaf over $V^*-D$. Here $D$ is the discriminant locus, consisting of $f$ such that $Y_f$ is singular. If $\Omega_f$ is globally defined on $V^*-D$, then its period integrals define a locally constant sheaf of finite rank over $V^*-D$.

Let $\pi: \cY\to V^*-D$ be the universal family of smooth CY hypersurfaces in $X$. For $f\in V^*-D$, let $j: Y_f\hookrightarrow X$ be the inclusion map. The vanishing cohomology of $Y_f$ is defined as $H^{d-1}(Y_f)_{van}:=Ker(H^{d-1}(Y_f)\xrightarrow{j_*}H^{d+1}(X))$. We have $H^{d-1}(Y_f)=j^*H^{d-1}(X)\oplus H^{d-1}(Y_f)_{van}$ and the splitting preserves the Hodge structures. We also have the following short exact sequence
$$
0\to H^{d}(X)_{prim}\to H^{d}(X-Y_f)\xrightarrow{Res}H^{d-1}(Y_f)_{van}\to 0.
$$

The groups $H^{d-1}(Y_f, \CC)$ form a flat vector bundle $\HH^{d-1}_{\CC}$ over $V^*-D$. Let $\cH^{d-1}=\HH_{\CC}^{d-1}\otimes_{\CC} \cO_{V^*-D}$ be the corresponding locally free sheaf over $V^*-D$. $\{H^{d-1}(Y_f)_{van}\}_{f\in V^*-D}$ form a local system over $V^*-D$ which we denote by $\HH^{d-1}_{van}$. Let $\cH^{d-1}_{van}=\HH^{d-1}_{van}\otimes_{\CC} \cO_{V^*-D}$. The map $f\mapsto Res~{\Omega_f}$ gives a section of $\cH^{d-1}_{van}$. Thanks to \cite{Griffiths}, the vanishing cohomology of a sufficiently ample hypersurface $Y\subset X$ can be realized as residues of meromorphic forms on $X$ with poles along $Y$. Moreover, this realization relates the Hodge level to the order of the pole.

\begin{thm}\cite{Griffiths}
Suppose $Y$ is a sufficiently ample hypersurface in $X$, i.e. for any $p\geq0$, $q>0$, and $s>0$ we have $H^q(X, \Omega^p_X(sY))=0$. The residues $Res~\frac{\eta}{f^k}$ for $\eta\in H^0(X, K_X(kY))$ generates $F^{d-k-1}H^{d-1}(Y)_{van}$. Here $F^{\bullet}H^{d-1}(Y)_{van}$ is the induced Hodge filtration on $H^{d-1}(Y)_{van}$.
\end{thm}

The result can be used to device a reduction procedure for computing differential equations for period integrals. However, the procedure is difficult to implement, except in simple examples. Inspired by mirror symmetry \cite{CDGP}, additional tools and alternative methods have been developed for hypersurfaces in a toric variety (see for e.g. \cite{Bat}\cite{HKTY1993} and references therein.)
In this case, one can explicitly construct a global family of meromorphic top forms $\Omega_f$, and a D-module that governs the period sheaf. The D-module turns out to be an extension of a GKZ hypergeometric system. General solutions  to GKZ systems and their holonomic ranks have been found \cite{GKZ1990}\cite{Adolphson}, under certain nondegeneracy conditions. In one important degenerate case (for applications in mirror symmetry) a closed formula for the general solutions near a particular singular point, a point of maximal unipotent monodromy, has also been constructed \cite{HLY1994}\cite{LLY2009}, giving explicit power series expansions for period integrals of CY hypersurfaces. We refer the reader to \cite{Saito}\cite{Cattani} for surveys and the extensive bibliography therein for studies on the GKZ hypergeometric systems in other important contexts.

The main motivation of the present paper is to study period integrals and deformations of CY complete intersections in a homogeneous space. In this paper, we shall mostly restrict ourselves to partial flag varieties. We begin, in section 2, by constructing period integrals for the universal family of these CY manifolds, by means of a global Poincar\'e residue formula. An explicit formula in the case of Grassmannians is given. We show that the family of CY is deformation complete. Next, we would like to explicitly construct, describe, and ultimately solve a D-module that governs the period integrals. One attempt would be to imitate the construction of GKZ systems for toric hypersurfaces. In this case, recall that the idea was to start with a ``natural'' basis of $H^0(X,K_X^{-1})$ (which is indexed by integral points of some polytope); relations of the integral points then give rise to GKZ type binomial differential operators that govern the period integrals. The torus action on $X$ yields additional first order operators. Together, the operators form a regular holonomic system. For homogeneous spaces, the standard monomial theory for representations of reductive groups provides a natural way to index bases of cohomology of line bundles over $X$. Thus one would expect that there should be a parallel approach to construct GKZ type systems in this case. Unfortunately, the D-modules one constructs this way are almost never holonomic -- there would not be enough binomial differential operators to determine the period integrals -- mainly because the variety defined by the binomial ideal typically has the wrong dimension in this case.

The present paper solves this problem by introducing a type of systems of differential equations, which we call {\it tautological systems}. For a fixed reductive algebraic group $G$, to every $G$-variety $X$ equipped with a very ample equivariant line bundle $L$ (or a list of such bundles), we attach a system of differential operators defined on $H^0(X,L)$, depending on a group character (section 3.) We show that the system is regular holonomic when $X$ is a homogeneous space. A number of examples, including a toric variety, are discussed (section 4.)  In this case, the construction recovers the GKZ hypergeometric system for CY hypersurfaces in a toric variety $X$, when $G$ is the usual torus acting on $X$. Likewise the extended version of GKZ system is recovered when $G$ is taken to be $Aut(X)$. Finally, we show that the period integrals of the universal family of CY complete intersections in a partial flag variety are solutions to a tautological system. We also give an explicit description of this system (section 5.) In section 6, we discuss some numerical examples and their solutions. Further generalizations and examples will appear in a future paper.

\noindent{\bf Acknowledgements.}
The first author would like to thank Peter Littelmann and Gerry Schwarz for helpful communications. We also thank the reviewer for helpful suggestions. B.H.L. is partially supported by NSF FRG grant DMS-0854965, and S.T.Y. by NSF FRG grant DMS-0804454.

\section{Poincar\'e residues for partial flag varieties}\label{residue}

We follow the standard convention that the Lie algebra of a group $H$ is denoted by the gothic letter $\fh$. 

\begin{thm}\label{residue}
Let $X=G/P$, where $G=SL_n$ and $P$ be a parabolic subgroup of $G$. Then there exists a nowhere vanishing, $G$-invariant, holomorphic form $\Omega$ of degree $d=\dim X$, defined on a principal $Z$-bundle over $X$ where $Z$ is an algebraic torus.
\end{thm}

\begin{proof}
We shall fix a Borel subgroup $B$ of $G$ and assume that $P\supset B$. Let $\Phi,\Delta=\{\alpha_1,..,\alpha_l\}$ ($n=l+1$) be respectively  the root system and the set of simple roots of $G$, relative to $B$. It is well-known that the  parabolic subgroups of $G$ containing $B$ are parameterized by subsets of $\Delta$. The set 
$$
S=\Delta-\{\alpha_{d_1},..,\alpha_{d_r}\},~~~(0=d_0<d_1<\cdots<d_r<d_{r+1}:=n)
$$ 
corresponds to the parabolic subgroup $P_S$, whose Lie algebra is
$$
\fp_S=\fb+\sum_{\alpha\in[S]}\fg_{-\alpha}
$$
Here $[S]$ is the set of positive roots in the linear span of $S$. In this case, the homogeneous space $X=G/P_S$ can be identified with the flag variety $F(d_1,..,d_r,n)$, which consists of $r$-step flags
$$
(0)\subset E_1\subset\cdots\subset E_r\subset\CC^n,~~~(\dim E_i=d_i.)
$$
Put
$$
\tilde M=M_{d_1,d_2}\times\cdots\times M_{d_r,n}
$$
where $M_{a,b}$ ($a<b$) denotes the space of $a\times b$ matrices of rank $a$. Put 
\begin{eqnarray}
K&=SL_{d_1}\times\cdots\times SL_{d_r}\cr
\hat K&=GL_{d_1}\times\cdots\times GL_{d_r}\cr
Z&=Z(\hat K)\cong(\CC^\times)^r.
\end{eqnarray} 
Let $(g_1,..,g_r,g)\in \hat K\times G$ act on $\tilde M$ by
$$
(g_1,..,g_r,g)\cdot(m_1,..,m_r)=(g_1m_1g_2^{-1},...,g_rm_r g^{-1}).
$$
The map $p:\tilde M\ra X$, $(m_1,..,m_r)\mapsto(R(m_1\cdots m_r)\subset\cdots \subset R(m_r))$, defines a $G$-equivariant principal $\hat K$-bundle over $X$. Here $R(m)$ denotes the span of the row vectors of the matrix $m$. The matrix entries of the $m_i$ are called the Stiefel coordinates of $X=G/P_S$.
Let $x_1,..,x_p$ ($p=\dim\hat K$) form a basis of holomorphic vector fields generated by the action of $\hat K$ on $\tilde M$. Let $\omega_0$ denote that coordinate top form on the matrix space $\tilde M$, regarded as an open subset of the affine space $\CC^{d_1d_2}\times\cdots\times\CC^{d_rd_{r+1}}$. The form $\omega_0$ is $G$-invariant and nowhere vanishing on $\tilde M$. Let $\iota_{x_i}$ be the interior multiplication operator with respect to the vector field $x_i$. Since the holomorphic vector fields $x_1,..,x_p$ are everywhere linearly independent, and since $\omega_0$ is nowhere vanishing, the holomorphic form
$$
\Omega=\iota_{x_1}\cdots\iota_{x_p}\omega_0
$$
of degree $\dim X$ is nowhere vanishing. Since $G$ and $K$ commutes, and since $\omega_0$ is $G$-invariant, $\Omega$ is also $G$-invariant. Since $K$ acts trivially on $\wedge^p\hat\fk$, and since $\omega_0$ is $K$-invariant, $\Omega$ is also $K$-basic (i.e. $K$-invariant and $K$-horizontal). Note that $\Omega$ is $Z$-horizontal, but not $Z$-invariant because $\omega_0$ is not $Z$-invariant. It follows that $\Omega$ defines a form on the $G$-equivariant principal $Z$-bundle $\tilde M/K\ra X$. This completes the construction of $\Omega$. 
\end{proof}

By construction, it is clear that $Z$ acts on $\CC\omega_0$, hence on $\Omega$, by the character $\prod_{i=1}^rdet(t_i)^{d_{i+1}-d_{i-1}}$. Note that
\begin{equation}\label{c1}
c_1(TX)=c_1(K_X^{-1})=\sum_{i=1}^r (d_{i+1}-d_{i-1})\lambda_{d_i}\in Pic(X)
\end{equation}
whose coefficients agrees with the exponents of the character above. Here $\lambda_1,..,\lambda_l$ are the fundamental dominant weights of $G=SL_n$.

\bs
\begin{exm}
$G(d,n)$, the Grassmannian of $d$-planes in $\CC^n$. 
\end{exm}
\cut

We will derive an explicit formula for $\Omega$ in this case. Let $\pi: M_{d,n}\to M_{d,n}/GL_d=G(d,n)$ be the natural projection. $\widehat{G}(d,n): =M_{d,n}/SL_d$ is a principal $\CC^*$-bundle over $G(d,n)$. Let $\{z_{ij}| 1\leq i\leq d; 1\leq j\leq n\}$ be Stiefel coordinates for $G(d,n)$. The $GL_d$-action on $M_{d,n}$ generates $d^2$ linearly independent vector fields $u_{ij}=\Sigma^n_{l=1}z_{il}\frac{\partial}{\partial z_{jl}} \quad (1\leq i,j \leq d).$ Let $\omega=\prod_{1\leq r\leq d; 1\leq s \leq n} d z_{rs}$. We will contract $\omega$ with $u_{ij}$ successively to obtain the desired $SL_d$-invariant $d(n-d)$-form $\Omega$.

Let $\iota_{u_j}=\iota_{u_{dj}}\cdots\iota_{u_{1j}}$ and $\iota_u=\iota_{u_d}\cdots  \iota_{u_1}$ for convenience. 

Let $I_{d,n}=\{I=(1\leq i_1<i_2<\cdots<i_d\leq n)\}$ and let $\det Z_I$ be the minor of $(z_{ij})_{k\times n}$ indexed by $I$, then
\begin{equation*}
\begin{array}{lll}
\iota_{u_j} & = & \iota_{u_{dj}}\cdots\iota_{u_{2j}}\iota_{u_{1j}}\\

 & = & \prod_{1\leq i\leq k}\iota_{(z_{i1}\frac{\partial}{\partial z_{j1}}+z_{i2}\frac{\partial}{\partial z_{j2}}+\cdots+z_{in}\frac{\partial}{\partial z_{jn}})}\\
 
 & = & \sum_{I\in I_{d,n}}\det Z_I\prod_{l\in I} \iota_{\frac{\partial}{\partial z_{jl}}}\\
\end{array} 
\end{equation*} 
\begin{equation*}
\begin{array}{lll}
\iota_u & = & (\sum_{I_d\in I_{d,n}}\det Z_{I_d}\prod_{l\in I_d} \iota_{\frac{\partial}{\partial z_{dl}}})\cdots (\sum_{I_1\in I_{d,n}}\det Z_{I_1}\prod_{l\in I_2} \iota_{\frac{\partial}{\partial z_{1l}}})\\

 & = & \sum_{I_j\in I_{d,n}}\det Z_{I_d}\cdots \det Z_{I_1}\prod_{1\leq j\leq d; l\in I_j} \iota_{\frac{\partial}{\partial z_{jl}}}\\
 
 & =& \sum_{\sigma_i\in S_n/S_d\times S_{n-d}}(-1)^{\sigma_1}\cdots(-1)^{\sigma_d}p_{\sigma_1(I_0)}\cdots p_{\sigma_d(I_0)}\prod_{1\leq j\leq d; l\in \sigma_j(I_0)} \iota_{\frac{\partial}{\partial z_{jl}}}
 \end{array} 
\end{equation*} 
where $I_0=(1<2<3<\cdots<d)$.

Applying the above contraction operator $\iota_u$ to $\omega$ gives us (up to a sign) 

$$\Omega= \sum_{\sigma_i\in S_n/S_d\times S_{n-d}}(-1)^{\sigma_1}\cdots(-1)^{\sigma_d}p_{\sigma_1(I_0)}\cdots p_{\sigma_d(I_0)}\prod_{1\leq r\leq d; d+1\leq s\leq n}dz_{r\sigma_r(s)} $$

 
Let $L_1,..,L_s$ be ample line bundles on the partial flag variety $X$ such that $L_1+\cdots+L_s=K_X^{-1}$. By the Borel-Weil theorem, the $H^0(X,L_i)$ are irreducible representations of $G$. 
Assume that $f_i\in H^0(X,L_i)$ are general sections.

\begin{cor}  The Poincar\'e residue $Res~{\Omega\over f_1\cdots f_s}$ defines a holomorphic top form on the CY manifold $Y$: $f_1=\cdots=f_s=0$.
\end{cor}
\begin{proof}
Recall that every line bundle $L=\sum n_{d_i}\lambda_{d_i}$ on $X$ is the tensor product of the pullbacks of line bundles $n_{d_i}\lambda_{d_i}$ on the Grassmannians $G(d_i,n)$, via the natural projections $X\ra G(d_i,n)$. For $L$ ample, by the Borel-Weil theorem, the sections of $L$ can be represented as polynomials of the Pl\"ucker coordinates of those Grassmannians. It is easy to check that the multi-degrees of these polynomials are given by $c_1(L)$, namely, $(n_1,..,n_r)$. Since the Pl\"ucker coordinates can be represented as polynomials of the Stiefel coordinates $(m_1,..,m_r)$, we can view the section $f_1\cdots f_s$ of $L_1+\cdots+ L_s$ as a $K$-invariant  polynomial of the Stiefel coordinates. Moreover, since $L_1+\cdots+L_s=K_X^{-1}$, this polynomial transforms under $\hat K$, by the same character as $\Omega$ does (cf. eqn. (\ref{c1}).) It follows that $\Omega\over f_1\cdots f_s$ is a $\hat K$-invariant meromorphic $d$-form on the matrix space $\tilde M$ (with $d=\dim X$.) Thus it defines a meromorphic top form on $X$. The pole of this form is clearly $\cup Y_i$, hence the form is a global section of $\Omega^t_X(\cup Y_i)$. This implies that $Res~{\Omega\over f_1\cdots f_s}$ defines a holomorphic top form on $Y$. 
\end{proof}

We briefly describe the cohomology of a CY complete intersection $Y$ in $X=F(d_1,..,d_r,n)$, with focus on the case of hypersurfaces for simplicity.  We begin with the cohomology of $X$ \cite{Hirzebruch}. The ring  $H^*(X, \ZZ)$ is  isomorphic to the quotient of 
$$S(x_1, x_2, \cdots, x_{d_1})\otimes S(x_{d_1+1}, x_{d_1+2}, \cdots, x_{d_2})\otimes\cdots\otimes S( x_{d_{r-1}+1}, \cdots, x_{d_r})$$ 
by the ideal generated by $S^+(x_1, x_2, \cdots, x_n)$, where $S(x_1, \cdots, x_k)$ is the ring of symmetric functions in $x_1, \cdots, x_k$ with integer coefficients, and $S^+$ is the set of elements of degree $>0$. 

Let $W_G, W_P$ be the Weyl groups of $G$ and $P=P_S$ respectively. The $B$-orbits on $G/P$ are in one-on-one correspondence with elements in the coset $W_G/W_P$. The Bruhat decomposition $G/P=\cup_{w\in W_G/W_P}BwP$ is a cellular decompositon and the cells have the form $BwP\cong \CC^{\mu(w)}$, where $\mu(w)$ is defined as follows.  Let $\Phi^+$ be the set of positive root of $G$. Let $\Psi=\Phi^+-[S]$ be the set of positive roots $\alpha$ such that $-\alpha$ is not a root for $P$. Then $\mu(w)$ be the number of elements in $w(\Phi^+)\cap (-\Psi)$. Thus the Poincar\'e polynomial of $X$ is $P_t(X)=\frac{1}{|W_P|}\sum_{w\in W_G} t^{2\mu (w)}$.

\bs
\begin{exm} \cite{GH}
CY hypersurfaces in Grassmannian $X=F(d,n)=G(d,n)$. 
\end{exm}
\cut

Fix a full flag $0=V_0\subsetneq V_1 \subsetneq V_2 \subsetneq \cdots \subsetneq V_n=\mathbb{C}^n$ in $\mathbb{C}^n$. For any sequence of integers $\underline{a}=(n-d\geq a_1 \geq a_2 \geq \cdots \geq a_d \geq 0)$, the associated Schubert cell is defined to be $$W_{a_1, a_2, \cdots, a_d}=\{\Lambda \in G(d,n) | \dim (\Lambda \cap V_{n-d+i-a_{i}})=i\}\cong \mathbb{C}^{d(n-d)-\sum a_i}$$ 
The closure $\overline{W}_{a_1, a_2, \cdots, a_d}=\{\Lambda \in G(d,n) | \dim (\Lambda \cap V_{n-d+i-a_{i}})\geq i\}$ is called the Schubert variety associated to $\underline{a}$ and has dimension $d(n-d)-\sum a_i$. The associated Schubert cycle is the homology class $\sigma_{a_1, a_2, \cdots, a_d}=[\overline{W}_{a_1, a_2, \cdots, a_d}]$.
The Schubert cells are all even dimensional cells and therefore all boundary maps are 0. This implies that the homology ring $H_{*}(X, \mathbb{Z})$ has no torsion and is freely generated by  all Schubert cycles $\sigma_{a_1, a_2, \cdots, a_d}$. In particular, the dimension of $H^{2p}(X,\ZZ)$ is equal to the number of Schubert cycles of codimension $p$, i.e. the number of sequences of integers $\underline{a}=(n-d\geq a_1 \geq a_2 \geq \cdots \geq a_d \geq 0)$ such that $\sum a_i=p$.

Let $Y$ be a CY hypersurface in $X=G(d,n)$. Its total Chern class is $c(TX)=\prod_{1\leq i\leq d; d+1\leq j\leq n}(1-x_i+x_j)$ and the total Chern class of the normal bundle is $c(N_{Y/X})=1+x_1+x_2+\cdots+x_d$. Thus one gets
\begin{equation}\label{e1}
\begin{array}{cll}
\chi(Y) & = & \int_X{c(TX)c_1(N_{Y/X})\over c(N_{Y/X})}\\
& =& \int_X \prod_{1\leq i\leq d; d+1\leq j\leq n}(1-x_i+x_j)\frac{(x_1+x_2+\cdots+x_d)}{1+x_1+x_2+\cdots+x_d}
\end{array}
\end{equation}

Since $Y$ is ample, by Lefschets hyperplane theorem and the fact that $H^{p,q}(X)=0$ for $p\neq q$, we have 
\begin{equation}\label{e2}
H^{p,q}(Y)\cong H^{p,q}(X)=0,~~~\mbox{ if}~~p\neq q~~\mbox{and}~~ p+q\neq \dim Y.
\end{equation}
Together with eqn. \ref{e1}, these conditions determine the Betti numbers of the hypersurface $Y$. The Hodge numbers of $Y$ can be computed as follows. Put $\wedge_y \Omega^1_X=\oplus_{p\geq 0}y^p \Omega^p_X$ and recall the $\chi_y$-genus: 
$$\chi_y(X)=\chi(X, \wedge_y \Omega^1_X)=\sum_{p,q} (-1)^q h^{p,q}(X)y^p.$$ By the Hirzebruch-Riemann-Roch formula, we have 
\begin{equation}
\begin{array}{cll}
\chi_y(Y) & = & \chi(X, \wedge_y \Omega^1_X(1+y[-Y])^{-1}(1-[-Y]))\\
& =& \int_X \frac{1-e^{-\lambda}}{1+ye^{-\lambda}}\prod_{1\leq i\leq d; d+1\leq j\leq n}\frac{(x_j-x_i)(1+ye^{x_i-x_j})}{1-e^{x_i-x_j}}
\end{array}
\end{equation}
where $\lambda=c_1([Y])$. Together with eqn.\ref{e2}, this equality determine the Hodge numbers of $Y$.

We now turn to the question of deformation completeness.
Let $L_1, L_2, \cdots, L_s$ be ample line bundles on $X$ such that $\sum_{i=1}^s L_i=K_X^{-1}$, as before. Put
$$V^*=H^0(X, \oplus_{i=1}^s L_i)=H^0(X, L_1)\times \cdots \times H^0(X, L_s)$$
and let $D\subset V^*$ be the locus of singular complete intersections. Consider  the universal family of smooth CY complete intersections $Y$ parameterized by $V^*-D$ in $X=F(d_1,..,d_r,n)$:
$$\cY=\{(p, f)\in X\times \PP (V^*-D)| f(p)=0\}\to \PP(V^*-D).$$  
Suppose $\dim Y>2$. Adapting an argument of \cite{bc}, we will show the deformation completeness and the discreteness of automorphisms of $Y$.

\begin{thm}
The family $\cY\to \PP(V^*-D)$ is a complete deformation. Moreover,
any fiber of this family has no nontrivial holomorphic vector fields. 
\end{thm}
\begin{proof}
$\forall f=(f_1, f_2, \cdots, f_s)\in V^*-D$, where $f_i\in H^0(X, L_i)$, let $Y$ be the complete intersection defined by $f=0$. Let $\kappa: T_f V^*\to H^1(Y, T_Y)$ be the Kodaira-Spencer map. 

The short exact sequence $0\to T_Y\to T_X|_Y \to \oplus_{i=1}^s\cO_Y(L_i)\to 0$ induces a long exact sequence
\begin{equation}\label{les3}
\begin{array}{rll}
0& \to & H^0(Y, T_Y)\to H^0(Y, T_X|_Y)\to H^0(Y, \oplus_{i=1}^s\cO_Y(L_i))\xrightarrow{\kappa} H^1(Y, T_Y)\\
& \to &  H^1(Y, T_X|_Y)\to \cdots
\end{array}
\end{equation}

To show that every small deformation of $Y$ is still a complete intersection, it suffices to show $\kappa$ is surjective, equivalently $H^1(Y, T_X|_Y)=0$, and that the restriction map $H^0(X, \oplus_{i=1}^s\cO_X(L_i))\to H^0(Y, \oplus_{i=1}^s\cO_Y(L_i))$ is surjective. We will need the following result \cite{bott}: 
\begin{equation}\label{bott}
H^q(X, \Omega^p_X)=0 \mbox{~if~} p\neq q \mbox{~and,~} H^p(X, T_X)=0 \mbox{~if~} p>0. 
\end{equation}

\begin{lem}
Let $E=\oplus_{i=1}^s L_i$ and $K^{-r}=\wedge^r E^*=\oplus_{1\leq j_1<\cdots< j_r\leq s}\cO(-L_{j_1}-\cdots -L_{j_r}))$. $K^0=\cO_X$. The Kozsul complex defined as follows is a resolution of $\cO_Y$.

$$K^\centerdot: 0\to K^{-s}\xrightarrow{d^{-s}} \cdots\xrightarrow{d^{-2}} K^{-1}\xrightarrow{d^{-1}} K^0$$

where $d^{-r}(e_1\wedge\cdots\wedge e_r)=\sum_{i=1}^r (-1)^{i-1}f_i e_1\wedge\cdots \wedge \widehat{e_i}\wedge \cdots\wedge e_r.$ 
\end{lem}

First, consider the two spectral sequences abutting to the hypercohomology of $K^\centerdot\otimes E$. 
We have $'A_1^{r,t}=H^t(X, K^r\otimes \cO_X(E))$ and $''A_2^{r,t}=H^t(X, \cH^r(K^\centerdot \otimes \cO_X(E)))$. By the above lemma, the second one degenerates and $$\mathbb{H}^0(X, K^\centerdot \otimes \cO_X(E))=''A_2^{1,0}=H^0(Y, \cO_Y(E)).$$ On the other hand, we will show that in the first spectral sequence $'A_1^{-r,r}=0$ for all $r\neq 0$, and thus $'A_1^{0,0}=H^0(X, \cO_X(E))$ maps onto $H^0(Y, \cO_Y(E))$. $'A_1^{-r,r}$ is different from 0 only for $0\leq r\leq s$ and 

\begin{equation}
\begin{array}{cll}
'A_1^{-r,r}& = & H^{r}(X, K^{-r}\otimes \cO_X(E))\\
& \cong& \oplus_{1\leq j_1<\cdots< j_r\leq s} H^{r}(X, \cO_X(-L_{j_1}-\cdots -L_{j_r})\otimes \oplus_{i=1}^s \cO(L_i))\\
& \cong& \oplus_{1\leq j_1<\cdots< j_r\leq s; 1\leq i\leq s} H^{r}(X, \cO_X(-L_{j_1}-\cdots -L_{j_r}+L_i)\\
& \cong& \oplus_{1\leq j_1<\cdots< j_r\leq s} H^{\dim X-r}(X, K_X\otimes\cO_X(L_{j_1}+\cdots +L_{j_r}-L_i))^*\\
\end{array}
\end{equation}

Since $K_X=-\sum_{j=1}^s L_j$, $K_X\otimes\cO_X(L_{j_1}+\cdots +L_{j_r}-L_i)$ is a negative line bundle. By the Kodaira vanishing theorem, the above cohomology group vanishes unless $r=0$. 

Second, consider the two spectral sequences abutting to the hypercohomology of $K^\centerdot\otimes T_X$. 
We have $'B_1^{r,t}=H^t(X, K^r\otimes T_X)$ and $''B_2^{r,t}=H^t(X, \cH^r(K^\centerdot \otimes T_X))$. The second one degenerates and $$\mathbb{H}^1(X, K^\centerdot \otimes T_X)=''B_2^{1,0}=H^1(Y, T_X|_Y).$$ On the other hand, we will show that in the first spectral sequence $'B_1^{-r,r+1}=0$ for all $r$, from which it follows that $H^1(Y, T_X|_Y)=0$.

\begin{equation}
\begin{array}{cll}
'B_1^{-r,r+1}& = & H^{r+1}(X, K^{-r}\otimes T_X)\\
& =& \oplus_{1\leq j_1<\cdots< j_r\leq s} H^{r+1}(X, \cO_X(-L_{j_1}-\cdots -L_{j_r})\otimes T_X)\\
& =& \oplus_{1\leq j_1<\cdots< j_r\leq s} H^{\dim X-r-1}(X, K_X\otimes\cO_X(L_{j_1}+\cdots +L_{j_r})\otimes \Omega^1_X)\\
\end{array}
\end{equation}
When $r\neq s$ and $r>0$, since $K_X\otimes\cO_X(L_{j_1}+\cdots +L_{j_r})$ is negative and thus $H^{\dim X-r-1}(X, K_X\otimes\cO_X(L_{j_1}+\cdots +L_{j_r})\otimes \Omega^1_X)=0$ by the Kodaira vanishing theorem. When $r=0$, $B_1^{0, 1}=H^{1}(X, T_X)=0$ by (\ref{bott}). When $r=s$, $B_1^{-s, s+1}=H^{\dim X-s-1}(X, \Omega^1_X)=0$ if $\dim X-s-1\neq 1$, i.e $\dim Y>2$. (This excludes the K3 case). 

Finally, since $K_Y$ is trivial, by Serre duality and by using the Lefschetz hyperplane theorem inductively we have $$H^0(Y, T_Y)\cong H^{\dim Y}(Y, \Omega^1_Y)^*\cong H^{1, \dim Y}(Y)^*\cong H^{1, \dim Y}(X)^*=0$$ 
and thus $Y$ has no nontrivial holomorphic vector fields. This completes the proof.
\end{proof}

Therefore $\mathbb{P}(V^*-D)//SL_n$ is a coarse moduli space for $Y$. It is unirational and $\dim\mathbb{P}(V^*-D)//SL_n=\dim V^*-n^2$ (\cite{dol1}, \cite{dol2}).

\section{Tautological systems}

In this section, let $G$ be a connected linear reductive algebraic group, let $X$ be a projective $G$-manifold. Let $V$ be a finite dimensional $G$-module, and let $\phi:X\into\PP V$ be a $G$-equivariant embedding. Denote by $I(X,\PP V)\subset\CC[V]=Sym~V^*$ the vanishing ideal of $\phi(X)$ in $\PP V$. Since we have a canonical symplectic form $\bra,\ket$ on $T^*V=V\times V^*$, each linear function $\zeta\in V^*$ uniquely defines a derivation $\partial_\zeta\in Der~\CC[V^*]$, by the formula $\partial_\zeta a=\bra a,\zeta\ket$, $a\in V$. The linear $G$ action $G\ra Aut~V$ induces a Lie algebra action $\fg\ra End~V\subset Der~\CC[V^*]$, $x\mapsto Z_x$, with $Z_x a=x\cdot a$. We refer to the $Z_x$ as the $G$ operators.

\begin{defn}\label{tauto}
 Fix an an integer $\beta$, viewed as a character of $\CC^\times$. The tautological system (or D-module) $\tau(X,\phi,G,\beta)$, is the $D_{V^*}$-module $D_{V^*}/J$ where $J$ is the left ideal of $D_{V^*}$ generated by the following operators: $\{p(\partial_\zeta)|p(\zeta)\in I(X,\PP V)\}$, $\{Z_x| x\in\fg\}$, together with the Euler operator $\sum_i a_i\partial_{\zeta_i}+\beta$, where the $a_i$ and $\zeta_i$ are dual bases of $V,V^*$. We shall call $I(X,\PP V)$ (and the ideal of differential operators corresponding to it,) the embedding ideal of $X$ relative to $\phi$.
\end{defn}

\begin{rem} The correspondence $\zeta\mapsto\partial_\zeta$ may be thought of as part of the Fourier transform between $D_{V}$ and $D_{V^*}$.
\end{rem}

The definition can be made slightly more general and purely algebraic by starting with the initial data: $V$ a $G$-module and an radical ideal $I\subset\CC[V]$, possibly inhomogeneous. Since only the action of the Lie algebra of $G\times\CC^\times$ on $V^*$ enters the definition, it can be extended to allow an arbitrary (hence possibly non-integrable) action of this Lie algebra, with the integer $\beta$ being replaced by a character. Since the main examples of this paper do come equipped with group actions, we will defer this generality to a future study. However, for application to complete intersections, we will later generalize the definition to allow multiple line bundles on $X$ (section 5.)

In Definition 1, p15 \cite{Hotta}, Hotta introduces the notion of an ``$L$-twistedly $G$-equivariant'' D-module  over a $G$-variety $X$, where $L$ is a connection on $G$. A result there can be applied to a tautological system.\footnote{We thank the reviewer for pointing this out to us.}. On p17 \cite{Hotta}, we can let $V^*$ play the role of $X$ there, $D_{V^*} I(X,\PP V)$ the role of $I$ there, and $\beta$ the role of the $G$-character $\lambda$. Then Hotta's Theorem 2 implies that a tautological system has the $\lambda$-twistedly $G$-equivariance property. In section 4, p22 of the same reference, he considers the special case when $X$ is given by the closure of the $G$-orbit of a point in a linear representation of $G$; this can also be viewed as a special case of Definition \ref{tauto}.


Using the dual bases, we can write an element of $ \CC[V]=\CC[\zeta_1,\zeta_2,...]$ as a polynomial $p(\zeta)=p(\zeta_1,\zeta_2,...)$, and $p(\partial_\zeta)$ as a partial differential operator $p({\partial\over\partial a_1},{\partial\over\partial a_2},...)$ with constant coefficients, acting on functions of the variables $a_1,a_2,...$. If $(x_{ji})$ is the matrix representing $x\in\fg$ acting on $V$ in the basis $a_i$, i.e. $x\cdot a_i=\sum_j x_{ji}a_j$, then 
$$Z_x=\sum x_{ji}a_j{\partial\over\partial a_i}.$$

Let $j:V\into W$ be a $G$-module homomorphism, $\pi:W^*\onto V^*$ the dual map. This induces a $G$-equivariant map on structure sheaves $\pi^\#:\cO_{V^*}\ra \pi_*\cO_{W^*}$, $f\mapsto f\circ \pi$ for $f\in\cO_{V^*}(U)$, and the induced homomorphism on germs is also a $\fg$-module homomorphism:
\begin{equation}\label{G-operator}
(Z_x^Vf)\circ\pi=Z^W_x(f\circ\pi)
\end{equation}
where $f\in\cO_{V^*}$.  Here $Z_x^V,Z_x^W$ are the $G$ operators on $V^*,W^*$ respectively. Likewise for the Euler operators (with the same character $\beta$.) 

Now $\pi:W^*\onto V^*$ induces the $G$-equivariant algebra homomorphisms $\pi:\CC[W]\onto\CC[V]$ and
$$
\pi:\C[\partial_{\zeta_W}|\zeta_W\in W^*]\onto
\C[\partial_{\zeta_V}|\zeta_V\in V^*]~~~\partial_{\zeta^W}\mapsto\partial_{\pi\zeta^W}.
$$
It is straightforward to check that for $f\in\cO_{V^*}$, $p(\partial_{\zeta^W})\in \C[\partial_{\zeta^W}]$, we have
\begin{equation}\label{poly-operator}
[(\pi p(\partial_{\zeta^W})) f]\circ\pi=p(\partial_{\zeta^W})(f\circ\pi).
\end{equation}

Let $\cS_{V^*}$ (likewise $\cS_{W^*}$) be the subsheaf of $\cO_{V^*}$ whose stalks consists of germs annihilated by the defining ideal $J$ of the D-module $\cM=\tau(X,\phi_V,G,\beta)$. Then we have a canonical isomorphism (of $\CC_{V^*}$-modules) from $\cS_{V^*}$ to the solution sheaf $\cH om_{\cD_{V^*}}(\cM,\cO_{V^*})$ of $\cM$. Under this identification, we can therefore view $\cS_{V^*}$ as the solution sheaf of $\cM$.

\begin{lem} \label{change-of-variables} (Change of variables)
Let $\phi_V:X\into\PP V$ be a $G$-equivariant embedding of a $G$-variety $X$, and $\phi_W=j\circ\phi_V$, where $j:V\into W$ is a $G$-module homomorphism, and $\pi:W^*\onto V^*$ its dual. Let $\cS_{V^*}\subset\cO_{V^*}$ and $\cS_{W^*}\subset\cO_{W^*}$ be the solutions sheaves of the D-modules $\tau(X,\phi_V,G,\beta)$ and $\tau(X,\phi_W,G,\beta)$. Then $\pi^\#$ maps $\cS_{V^*}$ isomorphically onto $\cS_{W^*}$.
\end{lem}
\begin{proof}
Since $\pi:\CC[W]/I(X,\PP W)\ra\CC[V]/I(X,\PP V)\cong\CC[X]$ is an isomorphism with $\pi I(X,\PP W)=I(X,\PP V)$,
it follows from eqns. (\ref{poly-operator}) and (\ref{G-operator}) that if $f\in\cS_{V^*}$, then $f\circ\pi\in\cS_{W^*}$.
In other words, $\pi^\#$ sends a solution to $\tau(X,\phi_V,G,\beta)$ to a solution to $\tau(X,\phi_W,G,\beta)$, so $\pi^\#:\cS_{V^*}\ra\cS_{W^*}$. It is injective, because it is so on the structure sheaves.

Fix bases $a_1,..,a_p$ of $V$, and $b_1,..,b_q$ of $W$, such that $j:a_i\mapsto b_i$ for $1\leq i\leq p$. Then we can regard sections of $\cO_{V^*}$, $\cO_{W^*}$ to be functions of the variables $a$ and $b$ respectively. Then $\pi^\#$ maps a function $f(a_1,...,a_p)$ on $V^*$ to $f(b_1,..,b_p)$, and $\pi$ maps ${\partial\over\partial b_i}$ to ${\partial\over\partial a_i}$ if $1\leq i\leq p$, and to zero otherwise. So, the ${\partial\over\partial b_i}$ ($p+1\leq i\leq q$) are generators in $\tau(X,\phi_W,G,\beta)$, hence they kills $\cS_{W^*}$. Thus all solutions on $W^*$ are independent of the variables $b_{p+1},..,b_q$.  Given any such solution $f(b_1,..,b_p)$, it is straightforward to check that the function $f(a_1,..,a_p)$ is a solution on $V^*$. This shows that $\pi^\#$ is surjective on the solution sheaves.
\end{proof}

The lemma will be used to give different descriptions to essentially the same D-module, by choosing different $G$-modules as targets for embedding $X$.
As the proof shows, the net effect of changing the target from $\PP V$ to $\PP W$ in the initial data of our tautological system, is that we introduce additional linear variables, and at the same time, additional first order operators corresponding to the linear forms in $V^\perp\subset W^*$.

Let $\phi:X\into\PP V$ be a given $G$-equivariant embedding. Let $\cM$ be the tautological $D_{V^*}$-module $\tau(X,\phi,G,\beta)$ for short. Let $\cH$ be the solution sheaf $\cH=\cH om_{D_{V^*}}(\cM, \cO_{V^*})$. The following is an analogue of Theorem 5.1.3 \cite{Kapranov1997}.

\begin{thm}\label{holonomicity}
Assume that the $G$-variety $X$ has only a finite number of $G$-orbits. Then the following statements hold.
\begin{enumerate}
\item The D-module $\cM$ is regular holonomic. In particular, there is an open subset $V^*_{gen}$ such that the restriction of the solution sheaf $\cH$ to $V^*_{gen}$ is locally constant of finite rank. 
\item More explicitly, let $X=\sqcup_{l=1}^rX_l$ be the decomposition into $G$-orbits and let $X_l^\vee\subset V^*$ be the conical variety whose projectivization $\PP(X_l^\vee)$ is the projective dual to the Zariski closure of $X_l$ in $X$. Then $V^*_{gen}=V^*-\cup_{l=1}^rX_l^\vee$.
\item Suppose the coordinate ring $\CC[X]$ is Cohen-Macaulay, then the rank of the solution sheaf $\cH$ over the generic stratum $V^*_{gen}$ is less than or equal to the degree of $X$ in $\PP V$. 
\end{enumerate}
\end{thm}
\begin{proof}
We will adopt a mix of arguments of Kapranov \cite{Kapranov1997} and Hotta \cite{Hotta}.

The Fourier transform of the tautological D-module $\cM=D_{V^*}/J$ is $\widehat{\cM}=D_V/\widehat{J}$, where $\widehat{J}$ is the $D_V$-ideal generated by $$I(X, \PP V), \{\sum x_{ji}\zeta_i{\partial\over\partial\zeta_j}+\sum x_{ii}, x\in\fg\}, \sum_i \zeta_i{\partial\over\partial\zeta_i}+\dim V-s$$
 It is a twisted $G\times \CC^\times$-equivarant coherent $D_V$-module in the sense of \cite{Hotta} whose support Supp $\widehat{\cM}$ is the cone over $X$ in $V$ and thus consists of finitely many $G\times\CC^\times$-orbits. Thus $\widehat{\cM}$ is regular holonomic \cite{Borel}.

The  tautological D-module $\cM=D_{V^*}/J$ is homogeneous since the ideal $J$ is generated by homogeneous elements under the graduation $\deg {\partial\over\partial a_i}=-1$ and $\deg a_i=1$. Thus $\cM$ is regular holonomic since its Fourier transform $\widehat{\cM}$ is regular holonomic (\cite{Brylinski}). 
\comment{The point here is that homogeneity guarantees that $\cM$ is monodromic in the sense of p82 \cite{SST}.}

The characteristic variety of $\cM$ has the following explicit description. Let $\widehat X_l$ be the cone of the $G$-orbit $X_l$ in $V-0$.  The group $\CC^\times$ acts on the cone $\widehat{X}$ over $\phi(X)$ by scaling, and $G\times\CC^\times$ acts on it with finite number of orbits $\widehat X_l$, together with the fixed point $\widehat X_0:=0$. Consider the characteristic variety $Ch(\cM)$ of $\cM$ as a subvariety of the symplectic variety $T^*V^*=V^*\times V$, equipped with the standard symplectic form $\bra,\ket$. $Ch(\cM)$ is contained in the following zero locus of principle symbols of the generators of $J$.
\begin{eqnarray*}
p(\zeta)&=&0,~~~\forall p\in I(X,V)\subset\CC[\zeta]\cr
\sum x_{ji}a_j\zeta_i&=&0,~~~x\in\fg\cr
\sum a_i\zeta_i&=&0.
\end{eqnarray*}
The first set of equations says that if $(a,\zeta)\in Ch(\cM)$, then $\zeta$ lies in $\widehat{X}$, hence in a unique orbit $\widehat{X}_l$. The second set of equations says that if $(a,\zeta)\in Ch(\cM)$ then $\bra a, Z_x^*\zeta\ket=0$, where $Z_x^*$ is viewed here as the tangent vector field corresponding to $x\in\fg$ generated by the dual $G$-action on $V$. The last equation says that $a$ is normal to the Euler vector field generated by the $\CC^\times$ action. In summary, if $(a,\zeta)\in Ch(\cM)$, then $\zeta$ lies in a $G\times\CC^\times$ orbit $\widehat X_l$ and $a\in V^*$ is normal to the orbit. In other words,
$$
Ch(\cM)\subset\sqcup_{l=0}^rT^*_{\widehat X_l} V
$$
where $T_{\widehat X_l}^* V$ is the conormal bundle of $\widehat X_l$ in $V$.

\begin{lem} By the natural identification $T^*V=T^* V^*$, we have
\begin{equation}
\sqcup_{l=0}^rT^*_{\widehat X_l} V=\cup_{l=0}^rT^*_{X_l^\vee} V^*\cup V\times\{0\}
\end{equation}\label{Whitney}
where $X_l^\vee\subset V^*$ be the conical variety whose projectivization $\PP(X_l^\vee)$ is the projective dual to the Zariski closure of $X_l$ in $X$. 
\end{lem}

\begin{proof}
The decomposition into $G\times \CC^\times $-orbits $\widehat{X}=\sqcup_{l=1}^r{\widehat X_l}\sqcup \{0\}$ is a Whitney stratification. In fact, for any pair $(\widehat{X_s}, \widehat {X_l})$ such that $\widehat {X}_s \subset \overline{\widehat{X}_l},$ there exists an open dense set $U\subset \widehat X_l$ such that $(\widehat X_s, U)$ satisfies Whitney condition; and since $G\times\CC^\times$-action on $U$ generates $\widehat X_l$, $(\widehat X_s, \widehat X_l)$ also satisfies Whitney condition. In particular, $\overline{T^*_{\widehat X_l} V}\cap \pi^{-1}(\widehat X_s)\subset T^*_{\widehat X_s} V$, where $\pi: T^*V\to V$ is the natural projection. Thus $$Ch(\cM)\subset\sqcup_{l=1}^rT^*_{\widehat X_l} V\sqcup \{0\}\times V=\cup_{l=1}^r\overline{T^*_{\widehat X_l} V}\cup \{0\}\times V.$$
$\overline{T^*_{\widehat X_l} V}\subset T^*V$ is a closed conical Lagrangian submanifold. Under the natural identification $T^*V^*=T^* V=V\times V^*$, we have $\overline{T^*_{\widehat X_l} V}=\overline{T^*_{X_l^\vee} V^*}$, where $X_l^\vee$ is the conical variety dual to $\overline{\widehat X_l}$, i.e. $\PP(X_l^\vee)\subset \PP V^*$ is the projective dual to $\overline{X_l}$.
\end{proof}
The singular locus of the D-module $\cM$ is the Zariski closure of the image of $Ch(\cM)-V^*\times \{0\}$ under the projection $T^*V^*\to V^*$ \cite{Saito} and it is contained in the union $\cup_{l=1}^rX_l^\vee$.  Thus the restriction of the solution sheaf $\cH$ to $V^*_{gen}=V^*-\cup_{l=1}^rX_l^\vee$ is a locally constant sheaf of finite rank. This proves parts (1)-(2).

To prove part (3), we apply the following lemma of Kapranov \cite{Kapranov1997}.
\begin{lem}\cite{Kapranov1997}
Let $\widehat{X}$ be the cone over $X$ in $V$, and let $E\subset V$ be a linear subspace such that $\dim E+\dim \widehat{X}=\dim V$ and $\widehat{X}\cap E=\{0\}$. If $\CC[X]$ is Cohen-Macaulay, then $$\dim_{\CC}(\CC[E]\otimes_{\CC[V]}\CC[X])=\deg X.$$ 
\end{lem}

The characteristic ideal $\widetilde{J}$ of $\cM=D_{V^*}/J$ is the ideal in $\CC[a, \zeta]$ generated by the principal symbols of differential operators in $J$. The holonomic rank of $\cM$ is $\rank(\cM)=\dim_{\CC(a)}(\CC(a)[\zeta]/\CC(a)[\zeta]\cdot\tilde{J}).$ This gives the rank of the solution sheaf $\cH$ over the generic stratum. 
\comment{See definition p31 [SST].}

For any point $a\in V^*_{gen}$, let $E_a$ be the linear subspace of $V$ defined by the linear equations $$\sum x_{ji}a_j\zeta_i=0,~~~x\in\fg;  \sum a_i\zeta_i=0.$$ Then $E_a\cap \widehat{X}=0$ by the construction of $V^*_{gen}$ and $\dim E_a+\dim \widehat{X}=\dim V$. Thus the quotient of $\CC[X]$ by the ideal generated by those linear equations has dimension equal to $\deg X$. This is the quotient obtained by moding out the principal symbols of generators of $J$ and  its dimension is bigger than or equal to the dimension of the quotient by the full characteristic ideal $\widetilde{J}$. Therefore $\rank(\cM)\leq \deg X$, i.e. the rank of the solution sheaf $\cH$ is less than or equal to $\deg X$.
\end{proof}

The homogeneous coordinate ring of the Grassmannian $G(d,n)$ is Cohen-Macaulay (see \cite{E}) and the degree of $G(d,n)\into\PP(\wedge^d\CC^n)$, under the Pl\"ucker embedding, is given by the following formula (see \cite{Harris}) $$
\deg G(d,n)=(d(n-d))!\prod_{i=0}^{d-1}\displaystyle\frac{i!}{(n-d+i)!}.$$ 
It follows that the degree of $\phi:X=G(d,n)\into\PP V$, $V=H^0(X,K_X^{-1})^*$, is this times $n^{d(n-d)}$.

Thus we have the following
\begin{cor}\label{Gkn-bound}
The tautological system $\tau(G(d,n),\phi,SL_n,1)$ is regular holonomic, and the solution sheaf $\cH$ is locally constant of finite rank over $V^*-X^\vee$, where $X^\vee\subset V^*$ is the discriminant locus parametrizing singular CY hypersurfaces in $G(d,n)$. Moreover, the rank of $\cH$ is less than or equal to $n^{d(n-d)}(d(n-d))!\prod_{i=0}^{d-1}\displaystyle\frac{i!}{(n-d+i)!}$.
\end{cor}

\section{Examples}

We keep the same notations as in preceding sections.

\begin{exm} Very ample equivariant line bundle.
\end{exm}
\cut

Let $X$ be a $G$-manifold, and $L$ be a very ample $G$-equivariant line bundle on $X$, and $\varphi_L:X\into\PP V$, $V=H^0(X,L)^*$, be the $G$-equivariant embedding provided by $L$. 
Note that the embedding ideal $I(X,\PP V)$ contains no nontrivial linear forms in this case. In particular, the Euler and the $G$ operators are the only first order generators of our tautological system in this case. The Change-of-variables lemma shows that if we can realize the same $G$-module $V$ inside another module $W$, then we obtain an alternative description for the solutions sheaf of $\tau(X,\varphi_L,G,\beta)$.

\begin{exm}\label{toric} Projective toric variety.
\end{exm}
\cut

Let $\cA=\{\bar\mu_0,...,\bar\mu_p\}\subset1\times\ZZ^n\subset\ZZ^{n+1}$ be a finite list of distinct vectors generating $\ZZ^{n+1}$, and $\cL\subset\ZZ^{p+1}$ be the lattice consisting of vectors $l$ such that $\sum_i l_i\bar\mu_i=0$. Note that $\sum_i l_i=0$. For each $l\in\cL$, put $l^\pm\in\ZZ^{p+1}_\geq$ such that $l=l^+-l^-$. Let $X_\cA$ be the projective variety defined by the  homogeneous (because $\sum_i l_i=0$) polynomials
$$
\zeta^{l^+}-\zeta^{l^-},~~(l\in\cL)
$$
in the variables $\zeta_0,..,\zeta_p$. This is an $n$ dimensional irreducible projective toric variety. The algebraic torus $T=(\CC^\times)^n$ acts on $X_{\cA}$ by
$$t\cdot [\zeta_0,..,\zeta_p]=[t^{\mu_0} \zeta_0,..,t^{\mu_p} \zeta_p]$$
where $\bar\mu_i=(1,\mu_i)$. 
Let $\phi:X_\cA\into\PP^p$ be the inclusion map, and $\beta$ any integer. The tautological D-module $\tau(X_\cA,\phi,T,\beta)$ coincides with a GKZ $\cA$-hypergeometric system, as introduced in \cite{GKZ1990}, and have important applications in mirror symmetry \cite{Bat}\cite{HKTY1993}. If we replace $T$ by $Aut~X_\cA$, the resulting tautological D-module becomes an extended GKZ $\cA$-hypergeometric system, as introduced in \cite{HLY1994}.

An important case that often arises in mirror symmetry is that one starts with a smooth projective toric variety $X$, such that $K_X^{-1}$ is semi-ample, i.e. $c_1(X)$ lies in the closure of the ample cone of $X$. The sections of the bundle defines a rational map $\varphi: X-\ra\PP V$ with $V=H^0(X,K_X^{-1})^*$, away from the base locus of the linear system. The closure of the image in $\PP V$ can then be identified with $X_\cA$ above, where $\cA$ can be explicitly determined.

\begin{exm}\label{Kapranov-systems} Kapranov's A-hypergeometric systems.
\end{exm}
\cut

In \cite{Kapranov1997} Kapranov introduced a generalization of the GKZ $\cA$-hyper -geometric systems by, roughly speaking, replacing the algebraic torus $T$ by a general reductive group $H$, a finite set of Laurent monomials by a finite set $A=\{V_\alpha\}$ of irreducible representations of $H$, and finally the ``index''  $\beta\in\ft^*$ by a character $\chi$ of the Lie algebra $\fh$. Put $M_A^*=\oplus_\alpha \End V_\alpha$ and let $\rho_A:H\ra\prod_\alpha GL(V_\alpha)\subset M_A^*-0$ be the direct sum module. Then the closure of $\rho_A(H)$ in $M_A^*-0$ is a spherical variety $Y_A$ on which $H\times H$ act naturally. Kapranov's A-hypergeometric system associated with the data $A$ is the differential system defined on the affine space $M_A$ given by (Eqn. 5.1.1 \cite{Kapranov1997}):
\begin{equation}\label{Kapranov A-system}
\left\{
\begin{array}{lll}
L_h\Phi & =  R_h\Phi =\chi(h)\Phi, & h\in \fh\\
P_f\Phi & = 0, & f\in I_A\\
\end{array} \right.
\end{equation}
Here $S^\bullet(M_A)\ra \CC[\partial]$, $f\mapsto P_f$, is the standard isomorphism between the polynomial ring on $M_A$ and the ring of differential operators on $M_A$ with constant coefficients;  $L_h$ and $R_h$ are the infinitesimal generators of the left and right actions of $H$ on $M_A$; $I_A$ is the ideal of $X_A:=\PP(Y_A)$. According to Definition \ref{tauto}, we see that Kapranov's system coincides with $\tau(X_A,\phi_A,H\times H,(-\chi,-\chi))$, where $\phi_A:X_A\into\PP M_A^*$ is the inclusion map. 

\begin{exm} D-modules with ``residual'' symmetry.
\end{exm}
\cut

Let $X$ be a $G$-variety, $V$ a $G$-module, and $\phi:X\into\PP V$ a $G$-equivariant embedding. Let $K$ be a closed subgroup of $G$. Then the D-module $\tau(X,\phi,K,\beta)$ admits a ``residual'' group action by the centralizer $H$ of $K$ in $G$. In particular, $H$ acts on the solution sheaf of this D-module. One interesting example is $X=\PP M$ where $M$ is the space of $n\times m$ ($n\geq m$) matrices, and $G=SL_n\times SL_m$ acting by the usual left-right multiplications. Put $K=SL_n\times 1$. Then there are exactly $m$ $K$-orbits in $X$ -- the matrices in $M$ of a given rank $(\geq 1)$ form a single orbit. In particular, our D-module is regular holonomic.
Since $X=\PP M$ is also a toric variety (under the action of a maximal torus of $SL(M)$), we can construct general solutions to $\tau(X,\cO(\dim M),SL(M),\beta)$ by using the method of \cite{HLY1994}. Since $K\subset G=K\times H\subset SL(M)$, our general solutions are \`a priori solutions to the D-module $\tau(X,\phi,K,\beta)$ which are also $H$-invariant.

\begin{exm}\label{partial-flag-variety} Partial flag variety.
\end{exm}
\cut

Let $X=F(d_1,..,d_r,n)$, and let $G=SL_n$ act on $\CC^n$ as usual. This induces a transitive $G$ action on $X$, and so there is exactly one $G$ orbit in $X$. Let $L$ be an ample line bundle on $X$, and put
$$
V=H^0(X,L)^*
$$
Then by the Borel-Weil-Bott theorem, $V$ is an irreducible representation of $G$. Moreover, given a highest vector $v\in V$, we have a unique $G$-equivariant embedding
$$
\varphi_L:X\into\PP V
$$  
which maps the standard flag in $X$ to $[v]$. By Theorem \ref{holonomicity}, it follows that the tautological system
$$
\tau(X,\varphi_L,G,\beta)
$$
is regular holonomic. 

In the next section, we show that for $L=K_X^{-1}$ and $\beta=1$, this system governs the period integrals of the universal family of CY hypersurfaces in $X$. Moreover, we will give an explicit description of this system by enumerating its generators.

\section{Tautological systems for partial flag varieties}

In this section, we shall study a tautological system associated to the partial flag variety 
$$X=F(d_1,..,d_r,n).$$
As in section 2, we fix a Borel subgroup $B$ of $G=SL_n$, and let 
$$
\Delta=\{\alpha_1,..,\alpha_{n-1}\}
$$
be the simple roots of $G$. Put
$$
S=\Delta-\{\alpha_{d_1},..,\alpha_{d_r}\}
$$  
and let $P_S$ be the parabolic subgroup of $G=SL_n$ corresponding to $S$.
Let $\Phi,\Phi^+$ be the corresponding the root system, and its positive roots.  Denote by $\lambda_1,..,\lambda_{n-1}$ the fundamental dominant weights of $SL_n$, so that $\bra\lambda_i,\alpha_j\ket=\delta_{ij}$. Thanks to result of \cite{bott}, the Picard group ($G$-equivariant or otherwise) of $X$ has the following description
 
\begin{prop} Put $X=F(d_1,..,d_r,n)$. Then
$Pic(X)=Pic_G(X)=Hom(P_S,\CC)\cong\sum_{i=1}^r\ZZ\lambda_{d_i}$. Moreover, a line bundle $L=\sum_{i=1}^r n_{d_i}\lambda_{d_i}$ is ample iff it is very ample iff $n_{d_i}\geq1$ for all $i$.
\end{prop}

By direct computation, we get

\begin{prop}\label{superample}
For $X=F(d_1,..,d_r,n)$, we have $K_X^{-1}=\sum_{i=1}^r n_{d_i}\lambda_{d_i}$ with $n_{d_i}\geq2$ for all $i$.
\end{prop}

Fix an ample line bundle $L=\sum n_{d_i}\lambda_{d_i}=\lambda\in Pic(X)$, and put
$$
V=H^0(X,L)^*.
$$
We can factor the embedding $\varphi_L:X\into\PP V$ in canonical way, as follows. 
We define the {\it incidence map}:
\begin{eqnarray*}
& &\iota: X=F(d_1,..,d_r,n)\into G(d_1,n)\times\cdots\times G(d_r,n)\cr
& & x\mapsto(\iota_1x,...,\iota_r x)
\end{eqnarray*}
where the $\iota_i:X\ra G(d_i,n)$ are the natural projection maps. For each $i$, we have the standard Pl\"ucker embedding of the Grassmannian 
$$\pi_i:G(d_i,n)\into\PP V_i$$ 
where $V_i=\wedge^{d_i}\CC^n$ is the fundamental representation of weight $\lambda_{d_i}$. 
We also have the Veronese maps 
$$\nu_i:\PP V_i\into\PP Sym^{n_{d_i}} V_i,~~[v]\mapsto[v\otimes\cdots\otimes v]$$
 and the Segre map 
$$
\psi:\PP Sym^{n_{d_1}}V_1\times\cdots\times\PP Sym^{n_{d_r}} V_r\into\PP W,~~([u_1],...,[u_r])\mapsto[u_1\otimes\cdots\otimes u_r]
$$
where 
$$W=Sym^{n_{d_1}}V_1\otimes\cdots\otimes Sym^{n_{d_r}} V_r.$$ 
Put 
$$\nu=\nu_1\times\cdots\times\nu_r,~~\pi=\pi_1\times\cdots\times\pi_r.$$
Then we get a $G$-equivariant embedding
$$
\phi=\psi\circ\nu\circ\pi\circ\iota:X\into\PP W
$$
such that $\phi^*O_{\PP W}(1)=L$. By the Borel-Weil theorem, $H^0(X,L)^*=V_\lambda$ is an irreducible module with highest weight $\lambda$.  Clearly, $\lambda$ is the highest weight in $W$ of multiplicity 1, implying that $W$ contains a unique copy of $V_\lambda$.
It follows the image $\phi(X)$ in $\PP W$ lies in the linear subspace defined by $V_\lambda^\perp\subset W^*$, consisting of the linear forms on $W$ annihilating $V_\lambda\subset W$. Moreover, $V_\lambda^\perp$ contains every linear form vanishing on $\phi(X)$. 

We now specialize to
$$
L=-K_X
$$
and proceed to describe the image $\phi(X)$ in $\PP W$. Put
$$
Y=\PP V_1\times\cdots\times\PP V_r.
$$
We first enumerates generators of the vanishing ideal $I(Y,\PP W)$ of $Y$ in $\PP W$, under the embedding 
$$\psi\circ\nu:Y\into\PP W.$$ 
Fix a basis $z_{ij}$ ($1\leq j\leq m_i=\dim V_i$) of $V_i^*=H^0(\PP V_i,O(1))$, and  introduce the notation
$$
z^v=z_1^{v_1}\cdots z_r^{v_r}=\prod_{i,j} z_{ij}^{v_{ij}},~~~v=(v_1,..,v_r)\in\ZZ_\geq^{m_1}\times\cdots\times\ZZ_\geq^{m_r}.
$$
Let $\cE$ be the set of exponent $v=(v_1,..,v_r)$ such that $|v_i|:=\sum_j v_{ij}=n_{d_i}$, for each component vector $v_i$. For $v\in\cE$, we can view $z^v$ as a monomial function on $V_1\times\cdots\times V_r$. Let $\xi_{i,v_i}$ ($|v_i|=n_{d_i}$) be the basis of  $H^0(\PP Sym^{n_{d_i}} V_i,O(1))$ such that
$\nu_i^*:\xi_{i,v_i}\mapsto z_i^{v_i}$ (the restriction map), and $\zeta_v$ the basis of $H^0(\PP W,O(1))$ such that $\psi^*:\zeta_v\mapsto\xi_{1,v_1}\cdots\xi_{r,v_r}$. 
It is easy to show that the binomial quadratic forms in $\PP W$
$$
\zeta_u\zeta_v-\zeta_w\zeta_t,~~~u+v=w+t,~(u,v,w,t\in\cE)
$$
vanish on $Y$. On the other hand it is also known that $I(Y,\PP W)$ is generated by quadratic forms. Then by term-wise elimination, we find that any quadratic form vanishing on $Y$ is a linear combination of the binomials above.
Finally, by Proposition \ref{superample}, we find that $I(X,\PP W)$ is generated by the linear forms $V_\lambda^\perp\subset W^*$, together with $I(Y,\PP W)$.

\begin{rem}
For general ample line bundle $L$, when the condition $n_{d_i}\geq2$ does not necessarily hold, the quadrics in $I(X,\PP W)$ can be much more complicated, involving quadratic forms which are not necessarily binomials. For example, if $X=G(d,n)$, $L=O(1)$ and $W=H^0(X,L)^*$, then $I(X,\PP W)$ is generated by the Pl\"ucker relations, which are of course not binomials.
\end{rem}

\begin{thm}\label{hypersurface}
Let $X=F(d_1,..,d_r,n)=G/P_S$. Let $K_X^{-1}=\sum_{\alpha\in\Delta-S} n_{d_i}\lambda_{d_i}$, and $\phi:X\into\PP W$, $W=Sym^{n_{d_1}}V_1\otimes\cdots\otimes Sym^{n_{d_r}} V_r$ be as defined above. Let $f\in H^0(\PP W,O(1))$ be a general section, and put $Y_f=\{f=0\}\cap\phi(X)$. For $\gamma\in H_{d-1}(Y_f,\ZZ)$, let $\tau(\gamma)$ be a tube over the cycle $\gamma$ in $X$. Then the period integral $\int_{\tau(\gamma)}{\Omega\over f}$ is a solution to the tautological system $\tau(X,\phi,G,1)$. The system is generated by the $G$-operators $Z_x$ ($x\in\fg$), the Euler operator $\sum a_v{\partial\over\partial a_v}+1$, the first order operators $\sum\bra a_v,\zeta\ket{\partial\over\partial a_v}$ ($\zeta\in V_\lambda^\perp\subset W^*$), together with the binomial operators
$$
{\partial\over\partial{a_u}}{\partial\over\partial{a_v}}-{\partial\over\partial {a_w}}{\partial\over\partial {a_t}},~~~u+v=w+t,~(u,v,w,t\in\cE.)
$$
\end{thm}
\begin{proof}
A general section has the form $f=f(a,\zeta)=\sum_v a_v\zeta_v$. By a direct calculation, for $p(\zeta)\in I(X,\PP W)$ of degree $s$, we find that
$$p(\partial_\zeta){1\over f}=(-1)^s s! {p(\zeta)\over f^{s+1}}$$
which is zero on $X$. This implies that the period integral is killed by $p(\partial_\zeta)$. Let $g$ be any automorphism of $X$. Since the period integral is the Poincar\'e pairing $\bra\tau(\gamma),{\Omega\over f}\ket$ on $X-Y_f$, it is invariant under $g$:
$$
\bra\tau(\gamma),{\Omega\over f}\ket=\bra (g_*)^{-1}\tau(\gamma),g^*{\Omega\over f}\ket.
$$
Now let $g\in G$ be close to identity. Then $(g_*)^{-1}\tau(\gamma)=\tau(\gamma)$. By Theorem \ref{residue}, $\Omega$ is $G$-invariant, so 
$$
\bra \tau(\gamma),{\Omega\over f}\ket=\bra\tau(\gamma),g^*({1\over f})\Omega\ket.
$$
For $x\in\fg$, consider the action of the 1-parameter subgroup $g=g_t=exp(tx)$ of $G$. We have
$${d\over dt}|_{t=0} g_t^*({1\over f})=-{x\cdot f\over f^2}=-{Z_xf\over f^2}=-Z_x({1\over f}).$$
 It follows that
$$
0=\bra\tau(\gamma),Z_x({1\over f})\Omega\ket=Z_x \bra \tau(\gamma),{\Omega\over f}\ket
$$
Finally, the period is killed by the Euler operator $\sum a_v{\partial\over\partial a_v}+1$ because ${1\over f}$ is homogeneous of degree $-1$ in the variables $a_v$. Thus we have shown that the period is killed by all generators of the tautological system.

The last assertion of the theorem follows from the argument preceding the statement of the theorem.
\end{proof}

By the Change-of-variables lemma, we can use the linear forms in $V_\lambda^\perp$ to eliminate variables by expressing the coordinate functions $a_v$ of $W^*$, in terms of a basis of $H^0(X,L)^*=V_\lambda$. Then the D-module $\tau(X,\phi,G,1)$ reduces to the D-module $\tau(X,\varphi_L,G,1)$, corresponding to the canonical embedding $\varphi_L:X\into\PP V_\lambda$. Thus, we can view the preceding theorem as giving an alternative description of $\tau(X,\varphi_L,G,1)$, by introducing more variables to the differential equation system, by factoring $\varphi_L$ in terms of the four classical maps, $\iota,\pi,\nu,\psi$. The reward is that the factorization gives a system whose quadratic operators are all {\it binomials} with simple (and universal) description, while the price is the introduction of a collection of first order operators with constant coefficients corresponding to the linear forms in $V_\lambda^\perp\subset W^*\subset I(X,\PP W)$.

Our results on hypersurfaces can be generalized to complete intersections as follows. In section 2, we have already seen the Poincar\'e residue formula for CY complete intersections in $X$. Let $L_i$ be $G$-equivariant ample line bundles, and $\beta_i\in\ZZ, i=1,2,\cdots, s$. Let $V_i=H^0(X, L_i)^*$ and $V=V_1\times\cdots\times V_s$. Let $a^i_l$ be a basis of $V_i$ and $\zeta^i_l$ be the dual basis of $V_i^*$. The line bundles $L_i$ defines an equivariant embedding $\phi:X\into \PP V_1\times\cdots\times\PP V_s$. Denote by $I(X,\PP V)$ the ideal of polynomial functions in $\CC[V]$ which vanish on $\phi(X)\subset\PP V_1\times\cdots\times\PP V_s$. Note that $I(X,\PP V_i)\subset I(X,\PP V)$. Let $Z^i_x, x\in \fg$ be the infinitesimal form of the $G$-action on $V_i$ and let $Z_x=\sum_{i=1}^sZ^i_x$.

\begin{defn}
We define the tautological system $\tau(X,G,L_1,..,L_s,\beta_1,..,\beta_s)$ as the $D_{V^*}$-module $D_{V^*}/J$, where $J$ is the left $D_{V^*}$-ideal generated by the following operators: $\{p(\partial_{\zeta})| p(\zeta)\in I(X,\PP V)\}$, $\{Z_x|x\in \fg\}$, and the Euler operators $\sum_la^i_l\partial_{\zeta^i_l}+\beta_i, i=1,2,\cdots, s$. 
\end{defn}

The argument for Theorem \ref{holonomicity} can be generalized to show that the system
$$
\tau(X,G,L_1,..,L_s,\beta_1,..,\beta_s)
$$ 
is a regular holonomic D-module. One considers the $G\times(\CC^\times)^s$ action on $V$, where $(\CC^\times)^s$ acts by scaling on the $s$ factors of the representation $V=V_1\times\cdots\times V_s$ of $G$. The characteristic variety of the D-module is then shown to be a subvariety of the disjoint union of conormal bundles of the finitely many $G\times(\CC^\times)^s$ orbits in the cone $\widehat{X}$ over $\phi(X)$, where each of the conormal bundles is Lagrangian in $T^*V^*$.

Now assuming $\sum_i L_i=K_X^{-1}$, generic sections $f_i\in H^0(X,L_i)$ define a smooth CY complete intersection $Y_f:=\{f_1=f_2=\cdots=f_s=0\}$ in $X$. For $\gamma\in H_N(Y_f,\ZZ)$ ($N=\dim Y_f$), let $\tau(\gamma)$  be a tube over the cycle $\gamma$ in $X$. Let $\beta_i=1, i=1,2,\cdots,s$. Then the period integral 
$$\int_\gamma Res{\Omega\over{f_1f_2\cdots f_s}}= \int_{\tau(\gamma)}{\Omega\over{f_1f_2\cdots f_s}}$$ 
is a solution to the tautological system $\tau(X,G,L_1,..,L_s,\beta_1,..,\beta_s)$. This follows from a verbatim argument as in the case of hypersurfaces (Theorem \ref{hypersurface}.)

\section{Explicit examples}


\begin{exm}{\it $X=\PP^{n-1}$.}
\end{exm}
\cut

Let $z_1, z_2, \cdots, z_n$ be homogeneous coordinates on $X$ and $$\Delta_1=conv \{ne_i\}_{i=1}^n\subset \mathbb{R}^n$$ be the set of exponents of degree $n$ monomials in $z_j$. $\Delta_1$ lies in an affine hyperplane in $\RR^n$. Shifting it by $(1,1,\cdots, 1)\in \ZZ^n$ and then projecting it to a coordinate plane, we get the following convex polytope in $\mathbb{R}^{n-1}$.

$$\Delta=conv \{(-1, \cdots, -1), (n-1, -1,\cdots, -1), \cdots, (-1, -1, \cdots, n-1)\}.$$

Let $A=\Delta\cap \ZZ^n=\{v_0,v_1, \cdots, v_{N-1}\}$ and suppose $v_0=(0,0,\cdots, 0)$. Let $t_j=\frac{z_j}{z_n},1\leq  j\leq n-1$ be affine coordinates on $X$. Then $$V^*=H^0(X, \cO(n))=\oplus_{i=1}^N \CC z_1z_2\cdots z_nt^{v_i}.$$ Let $\phi: X\hookrightarrow \PP V$ be the degree $n$ Veronese embedding. Let $$\cL=\{l\in \ZZ^N| \sum l_iv_i=0\}$$ be the lattice of integral relations among $v_i$. 

In this case, our construction in \ref{residue} recovers the well-known form
$$\Omega=\sum_{i=1}^n (-1)^{i-1} z_i dz_1 dz_2\cdots \widehat{dz_i}\cdots dz_n=z_1z_2\cdots z_n\prod_{i=1}^{n-1}\frac{dt_i}{t_i}.$$
Let $f(a, t)=z_1z_2\cdots z_n\sum a_i t^{v_i}\in V^*$ be a generic section such that $Y=\{f=0\}$ is smooth. Then $Res~\frac{\Omega}{f}$ is a nowhere zero holomorphic top form on $Y$. We have an isomorphism $H^{n-1}(X-Y)\xrightarrow[\cong]{Res} H^{n-2}(Y)_{van}$. When $n$ is odd, $H^{n-2}(Y)_{van}=H^{n-2}(Y)$.

Period integrals $\Pi_{\gamma} (a)= \int_{\tau(\gamma)} \frac{\Omega}{f(a,t)},$ where $\gamma \in H_{n-1}(X-Y)$ are solutions to the system $\tau(X, \phi, SL_n, 1)$, which coincides with the following extended GKZ-system \cite{HLY1994}:

\begin{equation}\label{extended-GKZ}
\left\{
\begin{array}{llll}
{\partial\over\partial a_i}{\partial\over\partial a_j}\Pi_{\gamma}(a) & = & {\partial\over\partial a_i'}{\partial\over\partial a_{j'}} \Pi_{\gamma}(a) &  if \ v_i+v_j=v_{i'}+v_{j'}\in \ZZ^n\\

Z_x \Pi_{\gamma}(a) & = & 0 & \forall x\in sl_n\\
\sum_i a_i{\partial\over\partial a_i} \Pi_{\gamma}(a) & =& -\Pi_{\gamma}(a) & \\
\end{array} \right.
\end{equation}

 Integrating $Res~\frac{\Omega}{f}$ along a particular cycle $$\gamma_0=\{|t_1|=|t_2|=\cdots=|t_{n-1}|=1\},$$ we get

\begin{equation}
\begin{array}{lll}
\Pi_{\gamma_0} (a)& = & \int_{|t_j|=1, \forall j} \frac{\sum_{j=1}^n (-1)^{j-1} z_j dz_1 dz_2\cdots \widehat{dz_j}\cdots dz_n}{\sum a_i z^{v_i}}\\
& = &  \int_{|t_j|=1,\forall j} \frac{1}{\sum a_i t^{v_i}}\prod_{i=1}^{n-1} \frac{dt_j }{t_j}\\
& =& \frac{1}{a_0} \sum_{l\in \cL, l_0<0, l_i\geq 0\ if\ i\neq 0} (-1)^{l_0} \frac{(-l_0)!}{\prod_{i\neq 0} l_i!} a^l
\end{array}
\end{equation}
\bigskip

Note that the period integral above gives a power series solution to the system \ref{extended-GKZ} near $a_0=\infty$, where $a_0$ is the coefficient of the monomial $z_1\cdots z_n$ in $f(a,z)$. By a result in \cite{HLY1996}, this is the only regular solution near this infinity. All other solutions have log singularities. An explicit formula for them can be given in terms of Gamma function. The formula also generalize to an arbitrary Fano toric manifold. See \cite{HLY1994}\cite{LLY2009} for details.

\begin{exm}{\it $X=G(2,4)$.}
\end{exm}
\cut

Under the Pl\"ucker embedding, $X=G(2,4)\hookrightarrow \PP^5$ is a quadratic hypersurface. 
The moduli space of CY hypersurfaces in $X$ is given by $\PP H^0(X, \cO(4))//SL_4$, whose dimension is 89.

A CY hypersurface in $X$ can be considered as a complete intersection of type $(2,4)$ in $\PP^5$ and the procedure to find periods on CY complete intersections in toric varieties applies. Let $z_0=p_{12}, z_1=p_{13}, z_2=p_{23}, z_3=p_{14}, z_4=p_{24}, z_5=p_{34}$ be homogeneous coordinates on $\PP^5$. Then $X=\{z_0z_5+z_2z_3-z_1z_4=0\}$. The particular period $$\Pi(a)=\int_{|z_0|=\cdots =|z_5|} \frac{\sum_i (-1)^i z_i dz_0\cdots \widehat{dz_i}dz_5}{(z_0z_5+z_2z_3-z_1z_4)(\sum a_i z^{v_i})},$$ where $\sum a_i z^{v_i}\in H^0(\PP^5, \cO(4))$, can be computed as a double residue. 

\smallskip

On the other hand, we can evaluate the period along a cycle in an affine chart in $G(2,4)$ explicitly as follows. The weight polytope of $SL_4$-action on $H^0(X, \cO(4))$ is $\Delta_w=conv\{4(e_i+e_j), 1\leq i<j\leq 4\}\subset \RR^4$. 

On the affine chart $U=\{( id_{2\times 2}, *_{2\times 2})\}\cong \CC^4$, we have $p_{12}=1$. $p_{13}=z_1, p_{23}=z_2, p_{14}=z_3, p_{24}=z_4$ are affine coordinates on this patch and $p_{34}=z_1z_4-z_2z_3$. We have a linear map $H^0(X, \cO(4))\to \oplus_{|v|\leq 8}\CC z^v$ expanding degree $4$ polynomials in $p_{ij}$ into polynomials in $z$ of degree $\leq 8$. Let $v_0$ be the exponent such that $p^{v_0}=p_{13}p_{23}p_{14}p_{24}=z_1z_2z_3z_4$.

The period integral $\Pi_{\gamma_0}(a)=\int_{\gamma_0}\frac{\Omega}{\sum a_i p^{v_i}},$ where $\sum a_i p^{v_i}\in H^0(X, \cO(4)),$ along the cycle $$\gamma_0=\{|z_{j}|=1,\ 1\leq j\leq 4\}$$ can be computed as follows.\footnote{We thank B. Song who first worked this out.}

\begin{equation}
\begin{array}{lll}
\Pi_{\gamma_0} (a)& = & \int_{|z_j|=1} \frac{z_1z_2z_3z_4}{\sum a_i p^{v_i}}\frac{dz_1dz_2dz_3dz_4}{z_1z_2z_3z_4}\\
& = &  \int_{|z_j|=1} \frac{1}{a_0+\sum_{i\neq 0} a_i p^{v_i-v_0}}\prod\frac{dz_j }{z_j}\\
& =& \frac{1}{a_0} \sum_{n=0}^\infty \textit{constant term in } (\frac{\sum a_i p^{v_i}}{z_1z_2z_3z_4})^n \frac{1}{a_0^n}\\
& =& \frac{1}{a_0} \sum_{n=0}^\infty \frac{n!}{\prod l_i!}a^l \cdot\textit{constant term in } \frac{p^{\sum l_i v_i}}{(z_1z_2z_3z_4)^n}\\
& =& \frac{1}{a_0} \sum_{n=0}^\infty \frac{n!}{\prod l_i!}a^l \cdot\textit{constant term in } \frac{p^{\sum l_i v_i}}{(z_1z_2z_3z_4)^n}\\
& =& \frac{1}{a_0} \sum_{l\in \cL, l_0\leq 0, l_i\geq 0 \ if\ i\neq 0} \frac{(-l_0)!}{\prod_{i\neq 0} l_i!} {{n_5(l)} \choose {n_5(l)+n_2(l)-n}} a^l
\end{array}
\end{equation}
where $n_5(l)$ is the exponent of $p_{34}$, $n_2(l)$ is the exponent of $p_{23}$, $n_1(l)$ is the exponent of $p_{13}$ in $p^{\sum l_i v_i}$ respectively. $\cL$ is the integral lattice defined by $\{\sum l_i w_i=0\}$, where $w_i$ is the weight of $(\CC^\times)^4$-action on $p^{v_i}$.

\begin{exm}{\it $X=G(d,n)$.}
\end{exm}
\cut

We have $Pic(X)=\mathbb{Z} \mathcal{O}(1)$ and $K_X^{-1}=\mathcal{O}(n)$, where $\cO(1)$ is the pullback of the hyperplane bundle via the Pl\"ucker embedding $X\hookrightarrow \PP(\wedge^d\CC^n)$. Any nonzero element in 
$V^*=H^0(X, \mathcal{O}(n))$ defines a CY hypersurface in $X$. Let $W^*=\oplus_{i=1}^N \CC p^{v_i}$ be the direct sum of all degree $n$ monomials in Pl\"ucker coordinates. Then $V^*=W^*/Q_n$, where $Q_n$ is the degree $n$ part of the ideal genereated by Pl\"ucker relations. Let $a_1, a_2, \cdots, a_N$ be coordinates on $W^*$.  

Consider the composition of the Pl\"ucker embedding and the Veronese embedding $\phi: X\hookrightarrow \PP(\wedge^d\CC^n)\hookrightarrow \PP W,$ which commutes with the anticanonical embedding $X\hookrightarrow \PP V\hookrightarrow \PP W$. The ideal defining $X$ in $\PP W$ is the sum of the ideals defining $\PP(\wedge^d\CC^n)$ and $\PP V$ in $\PP W$, and hence is generated by quadratic Veronese relations and linear relations in $Q_n$. Corresponding to these generators, we define the following differential operators. 
\begin{itemize}
\item
Veronese operators: $\{{\partial\over\partial a_i}{\partial\over\partial a_j}-{\partial\over\partial a_{i'}}{\partial\over\partial a_{j'}}| v_i+v_j=v_{i'}+v_{j'}\in \ZZ^{n \choose d}\}.$
\item degree 1 polynomial operators: $\{\sum_i c_i{\partial\over\partial a_i} |\sum_i c_ip^{v_i}\in Q_n\}$. 
\end{itemize}

Let $\rho: sl_n\to End(W^*)$ be the infinitesimal form of the $SL_n$-action on $W^*$. The vector field generated by any $x\in sl_n$ on $W^*$ can be written as

$$Z_x  =  \sum_{i,j} \rho_{ij}(x) a_j \frac{\partial}{\partial a_i}.$$
 
The period integrals on CY hypersurfaces in $G(d,n)$ are solutions to the following system of differential equations, which is an equivalent form of the tautological system $\tau(X, \phi, SL_n, 1)$, where $\phi: X\hookrightarrow \PP V$ is as above.
\begin{equation}
\left\{
\begin{array}{llll}
{\partial\over\partial a_i}{\partial\over\partial a_j}\Pi_{\gamma}(a) & = & {\partial\over\partial a_{i'}}{\partial\over\partial a_{j'}} \Pi_{\gamma}(a) &  if \ v_i+v_j=v_{i'}+v_{j'}\in \ZZ^{n \choose d}\\

\sum_i c_i{\partial\over\partial a_i}\Pi_{\gamma}(a) & = & 0 & \forall \sum_i c_ip^{v_i}\in Q_n\\

Z_x \Pi_{\gamma}(a) & = & 0 & \forall x\in sl_n\\
\sum_i a_i{\partial\over\partial a_i}\Pi_{\gamma}(a) & =& -\Pi_{\gamma}(a) & \\
\end{array} \right.
\end{equation}

\section{Note added: on holonomic rank}

Consider the universal family of CY hypersurfaces in a given partial flag variety $X$, Example \ref{partial-flag-variety}. The well-known applications of variation of Hodge structures in mirror symmetry show that it is important to decide which solutions of a differential system come from period integrals. The central object of our study is the period sheaf, i.e. the sheaf generated by the period integrals of the CY hypersurfaces. By Lemma \ref{change-of-variables} and Theorem \ref{hypersurface}, the period sheaf is a subsheaf of the solution sheaf of the module
$$\cM=\tau(X,\varphi_{K_X^{-1}}, SL_n,1).$$
Thus an important open problem is to decide when the two sheaves coincide. If they do not coincide, how much larger is the solution sheaf? From Hodge theory, it is well-known that the rank of the period sheaf is given by the dimension of the varying middle cohomology of the smooth hypersurfaces $Y_f$:
$$
\dim H^N(Y_f,\CC)-\dim i^*H^N(X,\CC)
$$
where $N=\dim Y_f$, and $i:Y_f\into X$ denotes the inclusion map. Therefore, to answer those questions, it is clearly desirable to know precisely the holonomic rank of $\cM$. 

Let us recall what is known about those general questions.
In the case of CY hypersurfaces in, say, a Fano toric manifold $X$, it is known \cite{GKZ1990}\cite{Adolphson} that the  rank of the GKZ hypergeometric system (cf. Example \ref{toric}) in this case is the normalized volume of the polytope generated by the exponents of the monomial sections in $H^0(X,K_X^{-1})$. This number is also the same as the degree of the anticanonical embedding $X\into\PP H^0(X,K_X^{-1})^*$. It is also known \cite{HLY1994} that this number always exceeds (and is usually a lot larger than) the  rank of the period sheaf. If one considers the extended GKZ hypergeometric system, where the torus $T$ acting on $X$ is replaced by the full automorphism group $\Aut X$, one would expect that  rank of this system should be closer to that of the period sheaf. In fact, based on numerical evidence, it was conjectured \cite{HLY1994} that for $X=\PP^n$ (which lives in both the toric world and the homogeneous world), the holonomic rank of $\cM$ coincides with the  rank of the period sheaf. In the case when $X=X_A$ is a spherical variety of a reductive group $G$ corresponding to a given set of irreducible $G$-modules $A$, Kapranov \cite{Kapranov1997} showed that the holonomic rank of his A-hypergeometric system (see Example \ref{Kapranov-systems}) is bounded above by the degree of embedding $X_A\subset\PP M_A^*$, if $\CC[Y_A]$ is assumed to be Cohen-Macaulay. Theorem \ref{holonomicity} generalizes this to an arbitrary tautological system $\tau(X,\phi,G,\beta)$ where $X$ has only a finite number of $G$-orbits.  Note, however, that the rank upper bound in each case cited above can be obtained without using (therefore does not take advantage of) assumptions about whether the underlying D-module arises from the variation of Hodge structures of CY hypersurfaces.

Since the release of the current paper in May 2011, progress has been made toward the problem of holonomic rank for CY hypersurfaces. We mention the following recent result.

\begin{thm}\cite{BHLY}
Let $X=G(d,n)$. Then the holonomic rank of the $\cD$-module $\cM$ at $f\in H^0(X,K_X^{-1})$ is precisely $\dim H^{d(n-d)}(X-Y_f)$.
\end{thm}

Note that the theorem holds not just at generic sections, but at every hyperplane section $f$.
The theorem has also been generalized to an arbitrary flag variety. The proof is beyond the scope of this paper, and is expected to appear shortly \cite{BHLY}. The theorem also implies the above mentioned conjecture of \cite{HLY1994}:

\begin{cor}
For $X=\PP^n$, the solution sheaf of $\cM$ coincides with the period sheaf of CY hypersurfaces in $X$.
\end{cor}

Let's compare this with the upper bound given by Corollary \ref{Gkn-bound}. In this case, $X=G(1,n+1)=\PP^n$ and the latter bound is $(n+1)^n$. We claim that this always exceeds the rank of the period sheaf, which is given by the dimension of the varying cohomology of a smooth CY hyperplane section in $X$. In fact, by using the Lefschetz hyperplane theorem, we find that the  rank of the period sheaf is
$$
\frac{n}{n+1}(n^n-(-1)^n)<n^n+1\leq(n+1)^n.
$$

\address{\SMALL B.H. Lian, Department of Mathematics, Brandeis University, Waltham MA 02454. lian@brandeis.edu.}

\address{\SMALL 
R. Song \& S-T. Yau, Department of Mathematics, Harvard University, Cambridge MA 02138. rsong@math.harvard.edu, yau@math.harvard.edu.}

\end{document}